\documentclass[twoside,reqno]{amsart}
\usepackage{pgfplots}
\usepackage{adjustbox}
\usepackage{amsfonts,amsmath,amscd,amsthm,amssymb}
\usepackage{mathtools}
\usepackage{url}
\usepackage[mathscr]{euscript} 
\usepackage{graphics}
\usepackage{wrapfig}
\usepackage{lscape}
\usepackage{rotating}
\usepackage{epsfig}
\usepackage{cite}
\usepackage{color}
\usepackage{booktabs}
\usepackage[vlined,ruled]{algorithm2e}
\usepackage{multirow}

\usepackage{tikz}
\usetikzlibrary{patterns}

\newtheorem{remark}{Remark}

\def\restrict#1{\raise-0.2ex\hbox{\ensuremath|}_{#1}}

\newcommand{\bld}[1]{\boldsymbol{#1}}

\newcommand{\Sh}{\bld{V}_{\!h}}

\newcommand{\divs}{{\nabla\cdot}}
\newcommand{\grads}{{\nabla}}

\newcommand{\pol}{\mathbb{P}}

\newcommand{\Oh}{{\mathcal{T}_h}}

\newcommand{\Eh}{\mathcal{F}_h}

\newcommand{\jmp}[1]{[\![#1 ]\!]}
\newcommand{\vertiii}[1]{{\left\vert\kern-0.25ex\left\vert\kern-0.25ex\left\vert #1
    \right\vert\kern-0.25ex\right\vert\kern-0.25ex\right\vert}}

\begin{document}
\title[Explicit divergence-free DG]{
An explicit divergence-free DG method for incompressible magnetohydrodynamics}
\author{Guosheng Fu}
\address{Division of Applied Mathematics, Brown University, 182 George St,
Providence RI 02912, USA.}
\email{Guosheng\_Fu@brown.edu}


\begin{abstract}
We extend the recently introduced explicit divergence-free DG scheme for incompressible hydrodynamics \cite{Fu18c}
to the incompressible magnetohydrodynamics (MHD).
A globally divergence-free finite element space is used for both 
the velocity and the magnetic field.
Highlights of the scheme includes global and local conservation properties, 
high-order accuracy, energy-stability, pressure-robustness.
When forward Euler time stepping is used,
we need 
two symmetric positive definite (SPD) hybrid-mixed Poisson solvers
(one for velocity and one 
for magnetic field) to advance the solution to the next time level.

Since we treat both 
viscosity in the momentum equation and resistivity in the magnetic induction equation 
explicitly, the method shall be 
best suited for inviscid or high-Reynolds number, low resistivity flows so that the CFL constraint is not too restrictive.
\end{abstract}
\maketitle

\section{Introduction}
\label{sec:intro}
The search for finite element methods that produce an exactly divergence-free velocity field 
for incompressible hydrodynamics has regained considerable interest in the past decade; see e.g. the recent review
article \cite{John17}. Equally well, the enforcement of the divergence-free constraint for 
the magnetic field in magnetohydrodynamics (MHD) is a long standing computational 
issue in numerical simulations, see e.g. \cite{Brackbill80,Toth,Dedner02}.

In this paper, we propose a divergence-free DG method for the incompressible MHD equation based on 
a {\it velocity-magnetic field} formulation, extending our previous work on a
divergence-free DG scheme for incompressible hydrodynamics \cite{Fu18c} to the incompressible MHD setting. 
In particular, we use a globally divergence-free finite element space 
for both the velocity and magnetic field. Hence, the 
pressure field and the divergence-free constraints in the equations are eliminated by design.
The scheme enjoys features such as global and local conservation properties, 
high-order accuracy, energy-stability, and pressure-robustness.
Moreover, the scheme can be efficiently implemented \cite{Fu18c} when coupled with {\it explicit} time stepping methods.
In particular, two hybrid-mixed Poisson solvers (equivalent to the mass matrix inversion of 
the divergence-free finite elements) is needed to advance solution in time when forward Euler time stepping is used.

We point out that while 
it is usually agreed that a Poisson solver for the evolution of velocity in incompressible flow 
is unavoidable, this might not be the case for the evolution of the magnetic field.
Indeed, when complete local finite element spaces were used for the 
magnetic field, e.g. the {\it locally divergence-free} DG method \cite{LiShu05}, 
explicit time stepping yields the inversion of a block-diagonal mass matrix,
hence, the Poisson solver for the magnetic field is avoided. 
However, such method does not provide an 
{\it exactly} divergence-free magnetic field due to the lack of $H(\mathrm{div})$-conformity, 
and reconstructing an exactly divergence-free magnetic field 
usually requires, again,  a Poisson solver \cite{Brackbill80}.
However, we particularly mention that the 
{\it constraint transport} (CT) method \cite{CT}
produces an exactly divergence-free magnetic field without the need of a global solver. 
We refer to the discussion in \cite{Toth} for a comparison of various version of CT methods.
See also the recent work on global divergence-free DG methods for 
compressible MHD \cite{Li12,Fu2018}, which can be interpreted as high-order 
CT-type methods. 
The key idea to achieve a divergence-free magnetic field for the CT methods
is to first advance in time the normal component of the magnetic field on the mesh interfaces, 
then apply a (locally defined) divergence-free reconstruction procedure.
It seems that such CT procedure is not available yet for the velocity evolution for incompressible flow.
For this reason, we chose to use the same 
Poisson-solver approach for the magnetic field evolution as that for the velocity evolution.
Our approach leads to a remarkably simple formulation which is provable energy stable in the semi-discrete 
case.
Moreover, its implementation can be trivially adapted from 
an existing explicit divergence-free DG code for incompressible hydrodynamics.

The rest of the paper is organized as follows.
In Section 2, the explicit divergence-free DG scheme is introduced for the 
incompressible inviscid MHD equations.
It is extended to the incompressible 
viscous MHD equations in Section 3.
Two dimensional numerical results are presented in Section 4.
Finally, we conclude in Section 5.

\section{Incompressible inviscid MHD}
\label{sec:model}
We consider the following incompressible inviscid MHD equations in conservative form:
\begin{subequations}
\label{euler}
 \begin{align}
  \label{euler-1}
\partial_t \bld u +\divs(\bld u\otimes \bld u -
\bld B\otimes \bld B)
+ \nabla (p+
\frac12
\bld B\cdot\bld B) = &\; 0, && \text{ in }\Omega,
  \\
  \label{euler-2}
\divs \bld u = &\; 0, && \text{ in }\Omega,
  \\
  \label{euler-3}
\partial_t \bld B +\divs(\bld u\otimes \bld B-\bld B\otimes \bld u) = &\; 0, 
&& \text{ in }\Omega.
 \end{align}
\end{subequations}
with initial condition 
\[
\bld u(x,0) = \bld u_0(x), 
\quad 
\bld B(x,0) = \bld B_0(x), 
\quad \forall x\in \Omega,
\]
where $\bld u$ is the velocity, $p$ is the pressure
, and $\bld B$ is the magnetic field, 
$\Omega\subset\mathbb{R}^d$(d=2,3) is a polygonal/polyhedral domain. 
The initial velocity and magnetic field  $\bld u_0(x)$ and $\bld B_0(x)$ are 
assumed to be divergence-free.
For simplicity, we consider periodic boundary conditions only.
However, the inflow/outflow/wall boundary conditions can be easily included, see \cite{Fu18c}.

We spefically mention that taking divergence of the equation \eqref{euler-3} yields 
\[
\partial_t(\divs \bld B) = 0.
\]
This implies that the $\divs \bld B=0$ condition will always be respected since 
the initial magnetic field satisfies $\divs \bld B_0=0$. 
The physical interpretation of this identity is that there exists no {\it magnetic monopoles}.
\subsection{Preliminaries}
Let $\Oh$ be a conforming simplicial triangulation of $\Omega$.
For any element $T \in\Oh$, we denote by $h_T$ its diameter and 
we denote by $h$ the maximum diameter over all mesh elements. 
Denote by $\Eh$ the set of facets of  $\Oh$.

\def\jump#1{[\![{#1}]\!]} 
\def\mean#1{\{\!\!\{{#1}\}\!\!\}} 

The following finite element space will be used to discretize both the velocity and 
magnetic field:
\begin{align}
\label{space}
\Sh^{k,-1} : =&\; \{\bld v\in \Pi_{T\in\Oh}[\pol^{k}(T)]^d, \;\;
\jmp{\bld v_T\cdot\bld n}_F = 0 \;\;\forall F\in\Eh,\quad
\divs \bld v=0 \text{ on } \Oh.\}
\end{align}
where the polynomial degree $k\ge 1$, 
and
$\jmp{\cdot}$ is the usual jump operator and $\pol^k$ the space of 
polynomials up to degree $k$. 
Notice that the space $\Sh^{k,-1}\subset H(\mathrm{div},\Omega)$ is not a standard finite element space in the sense 
that no {\it local} bases exist due to the divergence-free constraint. We introduce proper
Lagrange multipliers 
for the implementation of our scheme using such space, see details in 
\cite[Section 2.4]{Fu18c}; see also Remark \ref{rk1} below.

Finally, we introduce the jump and average notation. 
Let $\bld \phi_h$ be any function in $\Sh^{k,-1}$.
On each facet $F\in\Eh$ shared by two elements $K^-$ and $K^+$, we denote 
$(\bld \phi_h)^\pm|_F =\left.(\bld\phi_h)\right|_{K^\pm}$, and use
\begin{align}
 \label{avg-jmp2d}
\jump {\bld \phi_h}|_F  = \bld\phi_h^+ - \bld\phi_h^-,\quad \quad
\mean {\bld \phi_h}|_F  = \frac12(\bld \phi_h^++ \bld \phi_h^-)
\end{align}
to denote the jump and the average of $\bld \phi_h$. Here
$K^{-}$ is chosen such that $F$ is an outflow facet based on the velocity field $\bld u_h$, i.e.
$\bld u_h\cdot \bld n^- \ge0$ where $\bld n^-$ is the normal direction of $K^-$ on $F$.
Note that the jump term only contains jump for the {\it tangential} component of the vector $\bld \phi_h$ due to 
normal continuity of the finite element space.

\subsection{Spatial discretization}
The divergence-free space $\Sh^{k,-1}$ shall be used for both the velocity and 
magnetic field.
With this space in use, the divergence-free constraint \eqref{euler-2} is
point-wisely satisfied by design, 
and the (total) pressure gradient term in \eqref{euler-1} do not enter into the weak formulation of the 
scheme.
Using standard upwinding DG discretizations for the four convective terms 
\begin{align}
 \label{four}
 \divs(\bld u\otimes \bld u), \;\;
 \divs(\bld B\otimes \bld B), \;\;
 \divs(\bld u\otimes \bld B), \;\;
 \divs(\bld B\otimes \bld u),
\end{align}
we arrive at the following semi-discrete scheme:
find $(\bld u_h(t),\bld B_h(t))\in \Sh^{k,-1}\times \Sh^{k,-1}$
such that
\begin{subequations}
  \label{scheme-euler}
 \begin{align}
(\partial_t \bld u_h, \bld v_h)_{\Oh} + \mathcal{C}_{uu}(\bld u_h; \bld u_h, \bld v_h) 
-\mathcal{C}_{bb}(\bld B_h; \bld B_h, \bld v_h)=& 0,
\quad \forall \bld v_h\in \Sh^{k,-1},\\
(\partial_t \bld B_h, \bld \phi_h)_{\Oh} + \mathcal{C}_{ub}(\bld u_h; \bld B_h, \bld \phi_h) 
-\mathcal{C}_{bu}(\bld B_h; \bld u_h, \bld \phi_h)=& 0,
\quad \forall \bld \phi_h\in \Sh^{k,-1}.
\end{align}
where $(\cdot,\cdot)_\Oh$ denotes the standard $L^2$-inner product,
and the four convective operators are given below 
 \begin{align*}
\mathcal{C}_{uu}(\bld u_h; \bld u_h, \bld v_h):=&
\sum_{T\in\Oh}\int_T -(\bld u_h\otimes\bld u_h):\grads \bld v_h\,\mathrm{dx}
  +\int_{\partial T}(\bld u_h\cdot\bld n)({\bld u}_h^{-}\cdot\bld v_h)\,\mathrm{ds} \\
\mathcal{C}_{bb}(\bld B_h; \bld B_h, \bld v_h):=&
\sum_{T\in\Oh}\int_T -(\bld B_h\otimes\bld B_h):\grads \bld v_h\,\mathrm{dx}
  +\int_{\partial T}(\bld B_h\cdot\bld n)(\widehat{\bld B}_h\cdot\bld v_h)\,\mathrm{ds} \\
\mathcal{C}_{ub}(\bld u_h; \bld B_h, \bld \phi_h):=&
\sum_{T\in\Oh}\int_T -(\bld u_h\otimes\bld B_h):\grads \bld \phi_h\,\mathrm{dx}
  +\int_{\partial T}(\bld u_h\cdot\bld n)({\bld B}_h^{-}\cdot\bld v_h)\,\mathrm{ds}\\
\mathcal{C}_{bu}(\bld B_h; \bld u_h, \bld \phi_h):=&
\sum_{T\in\Oh}\int_T -(\bld B_h\otimes\bld u_h):\grads \bld \phi_h\,\mathrm{dx}
  +\int_{\partial T}(\bld B_h\cdot\bld n)(\widehat{\bld u}_h\cdot\bld \phi_h)\,\mathrm{ds},
 \end{align*}
where the four upwinding numerical fluxes are given as follows:
\end{subequations}
\begin{alignat}{2}
 \bld u_h^-|_F = &\bld u_h|_{K^-},\quad 
&& \quad \bld B_h^-|_F = \bld B_h|_{K^-},\\
 \widehat{\bld u}_h|_F =& \mean{\bld u_h}+\frac12s_F \jump{\bld B_h},&&\quad
 \widehat{\bld B}_h|_F = \mean{\bld B_h}+\frac12s_F \jump{\bld u_h},
\end{alignat} 
where $s_F$ takes value $1$ if $(\bld B_h\cdot \bld n)(\bld u_h\cdot\bld n)|_F>0$, and 
$-1$ if $(\bld B_h\cdot \bld n)(\bld u_h\cdot\bld n)|_F\le0$.
Notice that the numerical fluxes $\bld u_h^-$ and $\bld B_h^-$ are the 
 upwinding fluxes based on the velocity direction $\bld u_h\cdot\bld n$, whilst
 the numerical fluxes $\widehat{\bld u}_h$ and $\widehat{\bld B}_h$ are the 
 upwinding fluxes based on the magnetic field direction $\bld B_h\cdot\bld n$.
A simple calculation yields
\begin{align*}
 \mathcal{C}_{uu}(\bld u_h; \bld u_h, \bld u_h)=
\frac12 \sum_{F\in\Eh}\int_{F}|\bld u_h\cdot\bld n|\jump{\bld u_h}^2\,\mathrm{ds}
& \ge 0,\\
 \mathcal{C}_{ub}(\bld u_h; \bld B_h, \bld B_h)=
\frac12 \sum_{F\in\Eh}\int_{F}|\bld u_h\cdot\bld n|\jump{\bld B_h}^2\,\mathrm{ds}&
 \ge 0,\\
 \mathcal{C}_{bb}(\bld B_h; \bld B_h, \bld u_h)
 +\mathcal{C}_{bu}(\bld B_h; \bld u_h, \bld B_h){\color{white}xxxxxxxxxxx}&\\
=-\frac12\sum_{F\in\Eh}\int_{F}|\bld B_h\cdot\bld n|
 \left[
 \jump{\bld u_h}^2+\jump{\bld B_h}^2
 \right]\,\mathrm{ds}
& \le 0.
\end{align*}
Hence, a standard energy argument implies that 
the scheme \eqref{scheme-euler} is energy-stable in the sense that the total energy do not increase:
\begin{align}
&\partial_t 
\Big((\bld u_h,\bld u_h)_\Oh+
(\bld B_h,\bld B_h)_\Oh\Big)\nonumber\\
&=\sum_{F\in\Eh} \int_{F}(|\bld u_h\cdot\bld n|+|\bld B_h\cdot\bld n|)
 \left(
 \jump{\bld u_h}^2+\jump{\bld B_h}^2
 \right)\,\mathrm{ds}
\le 0.
\end{align}

\subsection{Temporal discretization}
The semi-discrete scheme \eqref{scheme-euler} can be written as
\begin{align}
\label{ode}
 \mathcal{M}(\partial_t\bld U_h) = \mathcal{L}(\bld U_h),
\end{align}
where $\mathcal{M}$ is the $2\times 2$ block-diagonal 
mass matrix for the compound space $\Sh^{k,-1}\times \Sh^{k,-1}$, and 
$\mathcal{L}(\cdot)$
the spatial discretization  operator, with $\bld U_h$ includes unknowns for both
$\bld u_h$ and $\bld B_h$.
Any explicit time stepping techniques can be applied to the ODE system  \eqref{ode}.
We use the following three-stage, third-order strong-stability 
preserving Runge-Kutta method (TVD-RK3) \cite{ShuOsher88} in 
our numerical experiments: 
\begin{align}
\label{rk3}
 \mathcal{M}\bld U_h^{(1)} = &\; \mathcal{M}\bld U_h^n +\Delta t^n \mathcal{L}(\bld U_h^n),\nonumber\\
\mathcal{M}\bld U_h^{(2)} = &\; \frac34\mathcal{M}\bld U_h^n +
\frac14\left[\mathcal{M}\bld U_h^{(1)}+\Delta t^n \mathcal{L}(\bld U_h^{(1)})\right],\\
\mathcal{M}\bld u_h^{n+1} = &\; \frac13\mathcal{M}\bld U_h^n 
+\frac23\left[\mathcal{M}\bld U_h^{(2)}+
\Delta t^n \mathcal{L}(\bld U_h^{(2)})\right],\nonumber
\end{align}
where $\bld U_h^{n}$ is the given solution at time level $t^n$ and 
$\bld U_h^{n+1}$ is the computed solution at time level $t^{n+1} = t^n+\Delta t^n$.
In each time step, $3\times 2=6$ inversion of the mass matrix for the divergence-free 
space $\Sh^{k,-1}$ are needed. 
Hence, the computational cost is essentially doubled per time step comparing with a 
corresponding hydrodynamic simulation.

\begin{remark}[Implementation]
\label{rk1}
The actual implementation of the scheme \eqref{rk3} that avoids global mass matrix inversion
was discussed in details in \cite{Fu18c} where appropriate Lagrange multipliers were introduced.
In particular, it was shown in
\cite[Section 2.4]{Fu18c} that the inversion of the mass matrix is 
equivalent to a (symmetric-positive-definite) hybrid-mixed Poisson solver.
\end{remark}

\section{Incompressible viscous MHD}
Now, we extend the scheme \eqref{scheme-euler} to 
the following incompressible, viscous, resistive MHD equations:
\begin{subequations}
\label{ns}
 \begin{align}
   \label{ns-1}
\partial_t \bld u +\divs(\bld u\otimes \bld u -
\bld B\otimes \bld B)
+ \nabla (p+
\frac12
\|\bld B\|^2) 
-\nu\triangle\bld u= &\; 0, && \text{ in }\Omega,
  \\
  \label{ns-2}
\divs \bld u = &\; 0, && \text{ in }\Omega,
  \\
  \label{ns-3}
\partial_t \bld B +\divs(\bld u\otimes \bld B-\bld B\otimes \bld u)
-
\eta
\triangle\bld B= &\; 0, 
&& \text{ in }\Omega,
 \end{align}
\end{subequations}
with divergence-free initial conditions 
\[
\bld u(x,0) = \bld u_0(x),\;\; \bld B(x,0) = \bld B_0(x)\quad \forall x\in \Omega,
\]
and periodic boundary conditions.
Here $\nu=1/\mathrm{Re}$ is the viscosity, and $\eta$ is the resistivity.

We discretize both viscous and resistive terms using a symmetric interior penalty DG method \cite{ArnoldBrezziCockburnMarini02}.
The semi-discrete scheme reads as follows:
find $\bld u_h(t)\in \Sh^{k,-1}$
such that
\begin{subequations}
 \begin{alignat}{2}
   \label{scheme-ns}
&(\partial_t \bld u_h, \bld v_h)_{\Oh} + \mathcal{C}_{uu}(\bld u_h; \bld u_h, \bld v_h) \nonumber\\
&\;\;\;\;-\mathcal{C}_{bb}(\bld B_h; \bld B_h, \bld v_h)
+\nu\mathcal{A}_{d}(\bld u_h, \bld v_h)&&= 0,
\quad \forall \bld v_h\in \Sh^{k,-1},\\
&(\partial_t \bld B_h, \bld \phi_h)_{\Oh} + \mathcal{C}_{ub}(\bld u_h; \bld B_h, \bld \phi_h) 
\nonumber\\
&\;\;\;\;-\mathcal{C}_{bu}(\bld B_h; \bld u_h, \bld \phi_h)
+
\eta
\mathcal{A}_{d}(\bld B_h, \bld \phi_h)
&&= 0,
\quad \forall \bld \phi_h\in \Sh^{k,-1}.
\end{alignat}
\end{subequations}
where the bilinear form for the diffusion term, $\mathcal{A}_{d}$, takes the following form
\begin{align}
  \label{viscous-1}
  \mathcal{A}_d(\bld u_h, \bld v_h)
  :=&\;\sum_{T\in\Oh}\int_T\nabla u:\nabla v\,\mathrm{dx}\nonumber\\
&  -\sum_{F\in\Eh}\int_F\mean{\nabla \bld u_h}\jump{\bld v_h\otimes\bld n}\,\mathrm{ds}\nonumber\\
&  -\sum_{F\in\Eh}\int_F\mean{\nabla\bld v_h}\jump{\bld u_h\otimes\bld n}\,\mathrm{ds}\nonumber\\
&  +\sum_{F\in\Eh}\int_F\frac{\alpha k^2}{h}\jump{\bld u_h\otimes\bld n}
  \jump{\bld v_h\otimes\bld n}\,\mathrm{ds},
\end{align}
with $\alpha>0$ a sufficiently large stabilization constant. We take $\alpha = 2$
in the numerics presented in the next section.

To obtain a fully discrete scheme, we use the same explicit stepping as the inviscid case. A standard 
CFL  time stepping restriction $\Delta t \le C\min\{h/v_{\max}, h^2/\max\{\nu,\eta\}\}$
is to be expected, where $v_{\max}$ is the maximum velocity speed. Hence, the method shall be applied to high Reynolds 
number, low resistivity flows where $\max\{\nu,\eta\} \ll 1$
to avoid severe time stepping restrictions.

Finally, we point out that when either $\nu$ or $\eta$ is not too small, 
one shall consider treating viscous/resistive terms implicitly to avoid parabolic time stepping. 
Hence, a Stokes solver is needed.
In that setting, an hybridizable discontinuous Galerkin (HDG) discretization \cite{CockburnGopalakrishnanLazarov09} of the 
diffusion operator is more favorable
in terms of the linear system solver efficiency; see some discussion on divergence-free HDG schemes 
for incompressible flow in \cite{Lehrenfeld:10,LehrenfeldSchoberl16}.

\section{Numerical results}
\label{sec:numerics}
In this section, we present some numerical tests in two dimensions to 
show the performance of our scheme.
The numerical results are performed using the NGSolve software \cite{Schoberl16}.
For the viscous operator \eqref{viscous-1}, 
we take the stabilization parameter $\alpha$ to be $2$. 
We use the TVD-RK3 time stepping \eqref{rk3} with sufficiently small time 
step size that ensure stability 
\[\Delta t \le \min\left\{ C\!F\!L_{\mathrm{conv}}h/v_{\max},C\!F\!L_{\mathrm{visc}} h^2/\max\{\nu,\eta\}\right\},\]
where $C\!F\!L_{\mathrm{conv}}$ and $C\!F\!L_{\mathrm{visc}}$ are 
the CFL stability constants depending on the polynomial
degree $k$.
We use a pre-factored sparse-Cholesky 
factorization for the hybrid-mixed Poisson solver that is needed in each time step.
\subsection*{Example 1: Accuracy test}
This example is used to check the accuracy of our schemes, 
both for the inviscid MHD equations \eqref{euler}
and for the viscous, ideal MHD equations 
\eqref{ns} with $\nu = 1/Re = 1/100$, $\eta=0$.
We take the domain to be $[0,2\pi]\times [0,2\pi]$ and 
use a periodic boundary condition.
The initial condition and source term are chosen such that 
the exact solution is  
\[
 u_1 =B_1= -\cos(x)\sin(y)\exp(-2t/\mathrm{Re}), 
 u_2 =B_2= \sin(x)\cos(y)\exp(-2t/\mathrm{Re}).
\]
The $L^2$-errors in velocity and magnetic field at $t=0.5$ on unstructured triangular meshes
are shown in Table \ref{table:error}.
It is clear to observe optimal $(k+1)$-th order of convergence for both variables in both cases.
\begin{table}[ht]
\caption{History of convergence of the $L^2$ errors
in velocity $\bld u_h$and magnetic field $\bld B_h$ at time 
$t=1$.} 
\centering 
\resizebox{1\columnwidth}{!}
{%
\begin{tabular}{ cc cc cc|cc cc}
\toprule
  & & \multicolumn{4}{c}{inviscid  MHD}
          & \multicolumn{4}{|c}{viscous, ideal MHD ($\mathrm{Re}=100$)}
\tabularnewline
$k$& $h$  
& $\bld u_h$-error & order & $\bld B_h$-error & order
& $\bld u_h$-error & order & $\bld B_h$-error & order
\tabularnewline
\midrule
\multirow{4}{2mm}{1} 
& 0.7854  &  1.853e-01  &  0.00 &  1.853e-01  &  0.00
   &  1.824e-01  &  -1.82 &  1.850e-01  &  4.57    \\ 
& 0.3927  &  4.119e-02  &  2.17 &  4.119e-02  &  2.17
   &  4.043e-02  &  2.17 &  4.111e-02  &  2.17    \\ 
& 0.1963  &  1.010e-02  &  2.03 &  1.010e-02  &  2.03
   &  9.892e-03  &  2.03 &  1.005e-02  &  2.03    \\ 
& 0.0982  &  2.463e-03  &  2.04 &  2.463e-03  &  2.04
   &  2.414e-03  &  2.03 &  2.444e-03  &  2.04    \\[1.5ex]
\multirow{4}{2mm}{2} 
& 0.7854  &  2.280e-02  &  0.00 &  2.280e-02  &  0.00
   &  2.191e-02  &  0.00 &  2.256e-02  &  0.00    \\ 
& 0.3927  &  2.783e-03  &  3.03 &  2.783e-03  &  3.03
   &  2.500e-03  &  3.13 &  2.622e-03  &  3.10    \\ 
& 0.1963  &  3.641e-04  &  2.93 &  3.641e-04  &  2.93
   &  2.920e-04  &  3.10 &  3.166e-04  &  3.05    \\ 
& 0.0982  &  4.554e-05  &  3.00 &  4.554e-05  &  3.00
   &  3.196e-05  &  3.19 &  3.698e-05  &  3.10    \\ [1.5ex] 
\multirow{4}{2mm}{3}
& 0.7854  &  1.223e-03  &  0.00 &  1.223e-03  &  0.00
   &  1.218e-03  &  0.00 &  1.239e-03  &  0.00    \\ 
& 0.3927  &  6.555e-05  &  4.22 &  6.555e-05  &  4.22
   &  6.888e-05  &  4.14 &  7.122e-05  &  4.12    \\ 
& 0.1963  &  4.242e-06  &  3.95 &  4.242e-06  &  3.95
   &  4.693e-06  &  3.88 &  4.955e-06  &  3.85    \\ 
& 0.0982  &  2.592e-07  &  4.03 &  2.592e-07  &  4.03
   &  2.802e-07  &  4.07 &  2.976e-07  &  4.06    \\ 
\bottomrule
\end{tabular}}
\label{table:error} 
\end{table}

\subsection*{Example 2: Orszag-Tang vortex problem}
We consider the Orszag-Tang vortex problem \cite{orszag1979}. 
The inviscid equation \eqref{euler}
on the domain $[0,2\pi]\times [0,2\pi]$ with a periodic boundary condition and an initial condition:
\begin{align}
 u_1(x,y,0) = -\sin(y)\quad & u_2(x,y,0) = \sin(x),\\
 B_1(x,y,0) = -\sin(y)\quad & B_2(x,y,0) = \sin(2x),
\end{align}
is solved.

\begin{figure}[ht!]
 \caption{\textbf{Example 2:} the computational mesh.
 Left: coarse mesh. Right: fine mesh.
 }
 \label{fig:ds}
 \includegraphics[width=.48\textwidth]{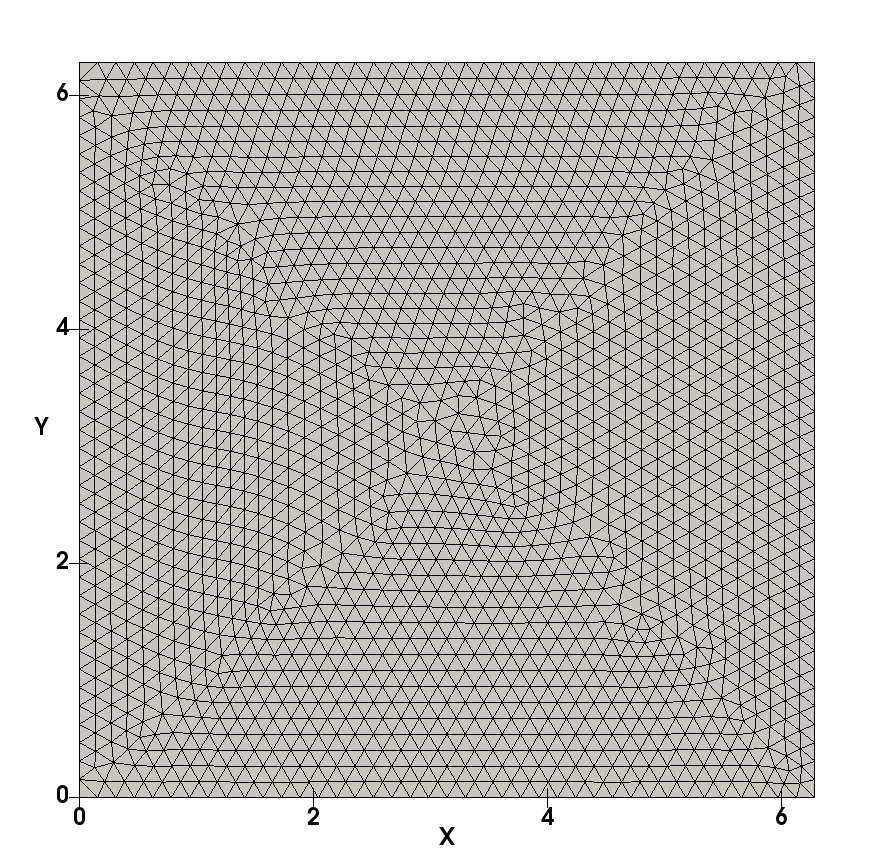}
 \includegraphics[width=.48\textwidth]{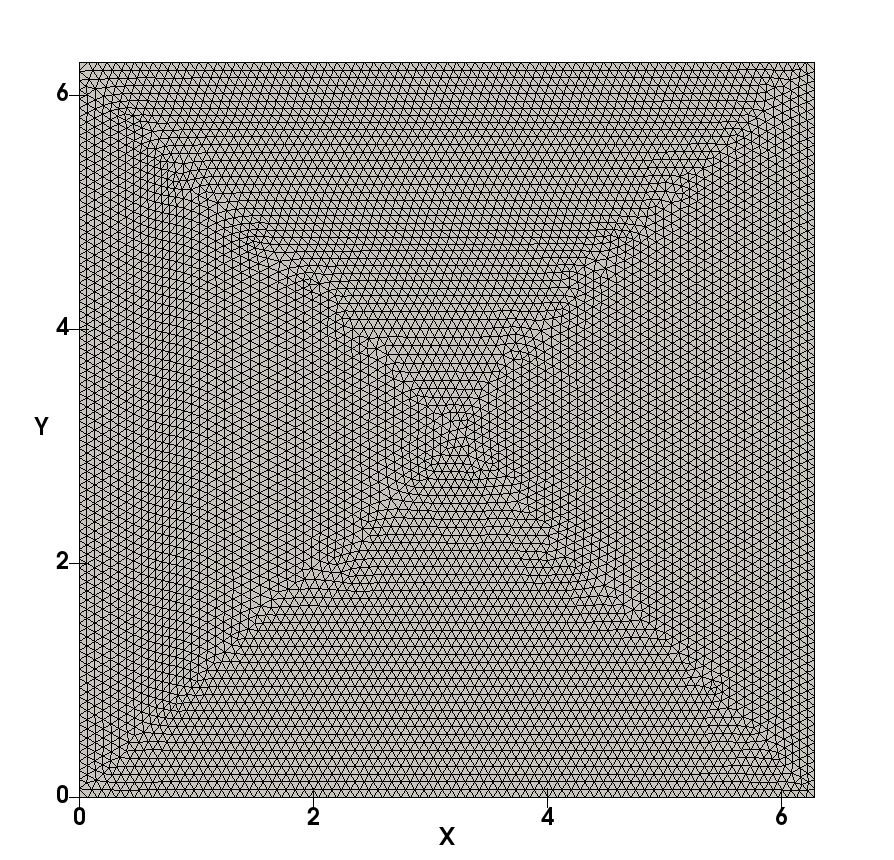}
\end{figure}

We use $P^3$ approximation on two
unstructured triangular meshes with mesh size $2\pi/40$ and $2\pi/80$, respectively,  
see Fig.~\ref{fig:ds}, and 
run the simulation  up to time $t=2$.

\begin{figure}[ht!]
 \caption{\textbf{Example 2:} Contour of vorticity.
 30 equally spaced contour lines between $-15$ and $15$.
 Left: results on the coarse mesh; right: results on the fine mesh.
 Top: $t=1$; bottom: $t=2$.
 $P^3$ approximation.
}
 \label{fig:ds1}
 \includegraphics[width=.48\textwidth, height=0.48\textwidth]{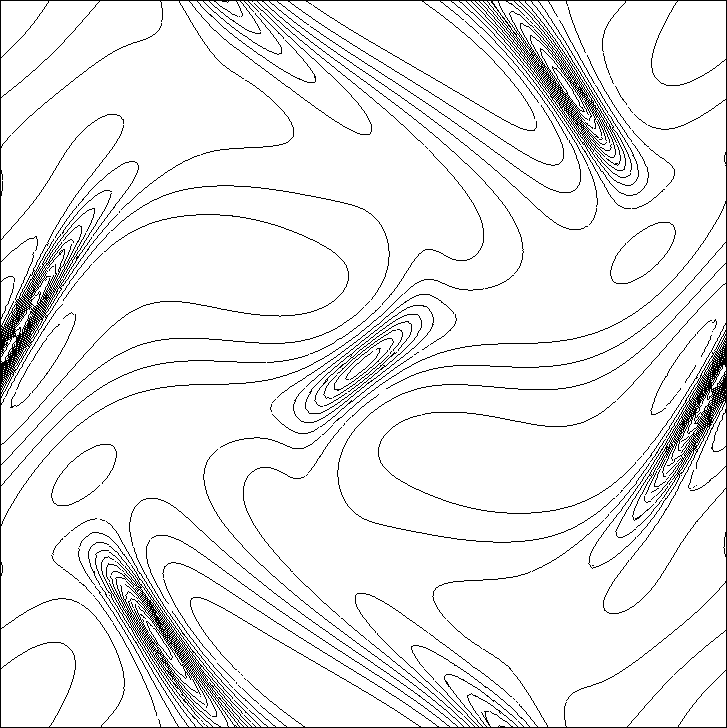}\;\;
 \includegraphics[width=.48\textwidth, height=0.48\textwidth]{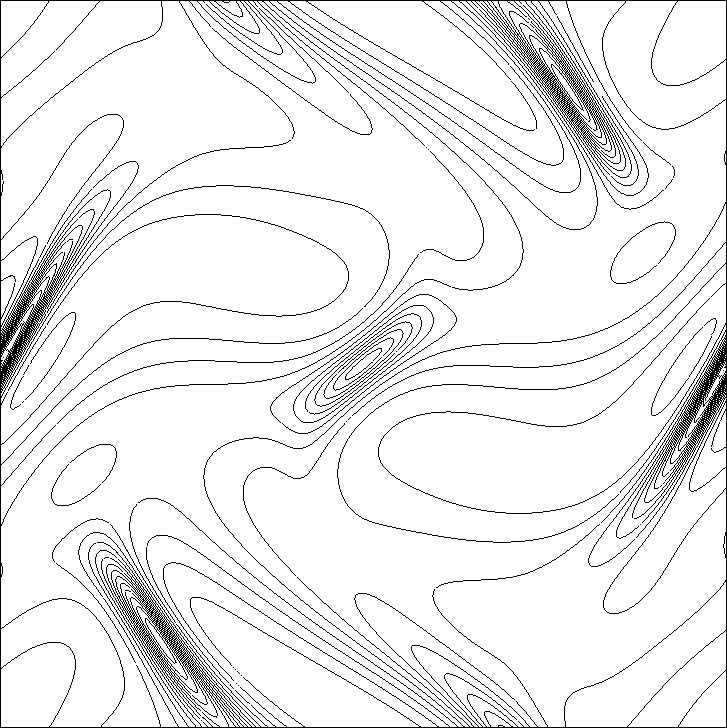}\\[1.5ex]
 \includegraphics[width=.48\textwidth, height=0.48\textwidth]{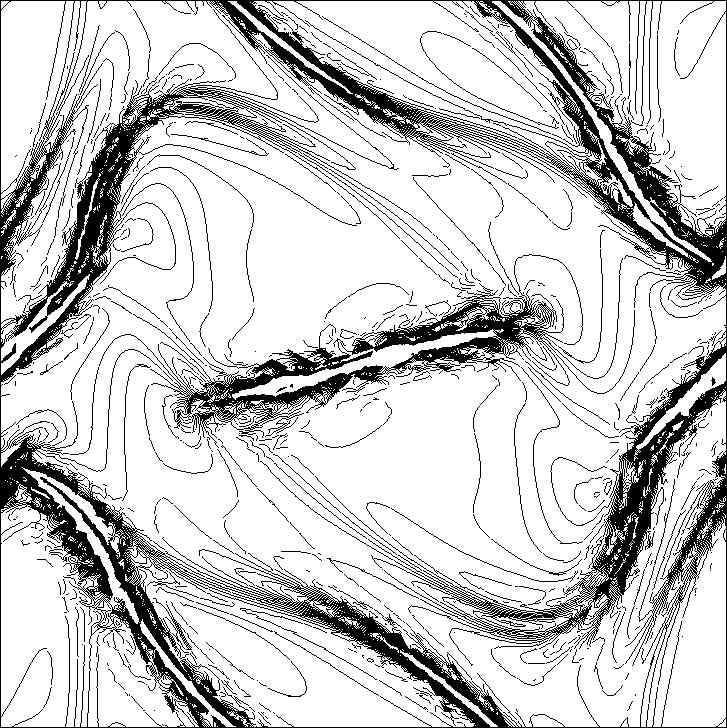}\;\;
 \includegraphics[width=.48\textwidth, height=0.48\textwidth]{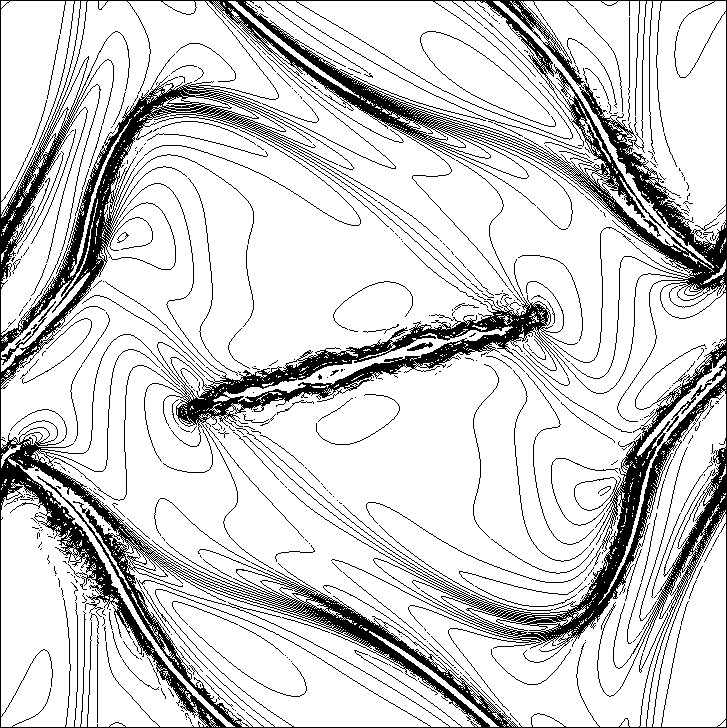}\\
\end{figure}

Contours of the vorticity $\omega_h:=\nabla_h\times\bld u_h$ at $t=1$ and $t=2$ are shown in  
Fig.~\ref{fig:ds1}. It is clear to observe the vorticity resolution improvement from 
the coarse mesh to the fine mesh at time $t=2$, where sharp unresolved layers have been developed.

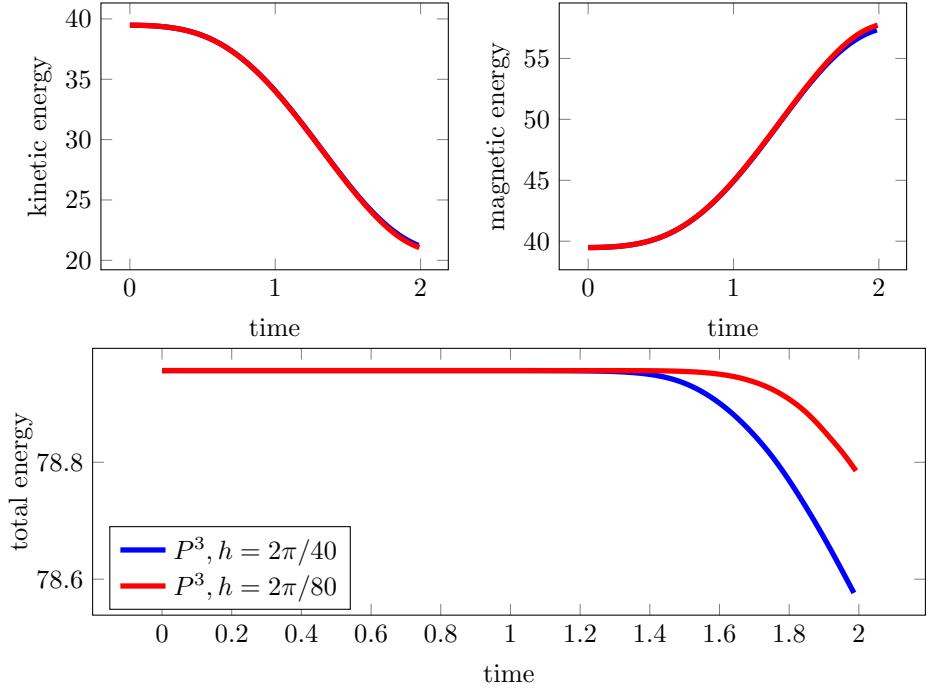
\begin{figure}[ht!]
 \caption{\textbf{Example 2:} 
the time history of kinetic, magnetic, and total energy. 
 }
 \vspace{1ex}
 \label{fig:ds0}
 \begin{tikzpicture} 
 \begin{axis}[
 	width=0.49\textwidth,
 	height=0.25\textheight,
 	xlabel={time},
    	ylabel={kinetic energy},
	      y label style={at={(axis description cs:0.15,.5)},anchor=south},
    	every axis plot/.append style={line width=2pt, smooth},
    	no markers,   	
    	legend style={at={(0.02,0.02)},anchor=south west}
 	]
 \addplot table[x index=0, y index=1]{otp3c.csv};
 \addplot table[x index=0, y index=1]{otp3f.csv};
 \end{axis}
 \end{tikzpicture}
 \hspace{1ex}
  \begin{tikzpicture} 
 \begin{axis}[
 	width=0.49\textwidth,
 	height=0.25\textheight,
 	xlabel={time},
    	ylabel={magnetic energy},
      y label style={at={(axis description cs:0.15,.5)},anchor=south},
    	every axis plot/.append style={line width=2pt, smooth},
    	no markers,   	
    	legend style={at={(0.02,0.02)},anchor=south west}
 	]
 \addplot table[x index=0, y index=2]{otp3c.csv};
 \addplot table[x index=0, y index=2]{otp3f.csv};
 \end{axis}
 \end{tikzpicture}
  \hspace{.5ex}
 \begin{tikzpicture} 
\begin{axis}[
	width=1\textwidth,
	height=0.25\textheight,
	xlabel={time},
	ylabel={total energy},
	      y label style={at={(axis description cs:0.05,.5)},anchor=south},
   	yticklabel style={/pgf/number format/fixed,/pgf/number format/precision=4},   	
   	every axis plot/.append style={line width=2pt, smooth},
   	no markers,   	
   	legend style={at={(0.02,0.02)},anchor=south west}
	]
\addplot table[x index=0, y index=3]{otp3c.csv};
\addplot table[x index=0, y index=3]{otp3f.csv};
\addlegendentry{$P^3, h = 2\pi/40$}
\addlegendentry{$P^3, h = 2\pi/80$}
\end{axis}
\end{tikzpicture}
\end{figure}

We then plot the time history of 
kinetic energy $(\bld u_h, \bld u_h)_\Oh$, 
magnetic energy $(\bld B_h, \bld B_h)_\Oh$, and the 
total energy (kinetic+magnetic) in Fig.~\ref{fig:ds0}.
We observe an energy transformation from kinetic energy to magnetic energy.
We also observe that the total energy is monotonically decreasing, which is to be expected since 
both our spatial and temporal discretization are dissipative.
The dissipated energy at time $t=2$ for the 
scheme on the coarse mesh is about $5\%$, 
while that on the fine mesh is about $2.5\%$.


\subsection*{Example 3: MHD Kevin-Helmholtz instability problem}
We consider an incompressible MHD Kevin-Helmholtz instability problem, the set-up is
adapted from \cite{Frank96} where a compressible MHD Kevin-Helmholtz instability problem was studied.
The inviscid equations \eqref{euler} on the domain $[0,L]\times [0,L]$ with $L=1$ 
is solved with a periodic boundary condition on 
the $x$-direction, and the slip wall boundary condition for both the velocity and magnetic field
$\bld u\cdot \bld n=\bld B\cdot \bld n=0$ in the  $y$-direction.
The initial conditions are taken to be 
\begin{align*}
 u_1(x,y,0) =&\; -\frac{u_0}{2}\, \mathrm{tanh}\left(\frac{y-L/2}{a}\right)
 +c_n\partial_y \psi(x,y),\\
 u_2(x,y,0) =&\; -c_n\partial_x \psi(x,y),\\
 B_1(x,y,0) =&\; \frac{u_0}{M_A},\\
 B_2(x,y,0) =&\; 0,
\end{align*}
with corresponding stream function
\[
 \psi(x,y) = u_0\exp\left(-\frac{(y-L/2)^2}{a^2}\right)\cos\left(\frac{2\pi}{L} x\right).
\]
Here, $a=L/20=0.05$ is the velocity shear scale length, 
$c_n = 10^{-3}$ is the noise scaling factor, and
$M_A$ is the Alfv\'enic Mach number.
Similar as in \cite{Frank96}, we take $M_A = 2.5$ (strong magnetic field) and $M_A = 5$ (weak magnetic field) 
in our numerical simulation. For comparison purpose, we also present numerical results for the 
hydrodynamic case (corresponding to $M_A=+\infty$).
Introducing the scaled time $\bar t = t/\Gamma$, with $\Gamma = 0.106{u_0}/{2a}=1.06$,
we run the simulation till scaled time $\bar t = 6$ (corresponding to physical time $t =6.36$).

For all the numerical tests in this example, we consider a uniform rectangular mesh with $256\times 256$ cells.
The divergence-free finite element space \eqref{space} on rectangular meshes
is modified to be 
\begin{align*}
\Sh^{k,-1} : =&\; \{\bld v\in 
\Pi_{T\in\Oh}\mathbf{RT}_k(T), \;\;
\jmp{\bld v_T\cdot\bld n}_F = 0 \;\;\forall F\in\Eh,\quad
\divs \bld v=0 \text{ on } \Oh\},
\end{align*}
where $\mathbf{RT}_k(T) := [P_{k+1,k}(T), P_{k,k+1}(T)]$ is the usual 
Raviart-Thomas \cite{RaviartThomas77} space on a rectangle $T$ with 
$P_{m,n}(T)$ the space of polynomials with degree at most $m$ in the $x$-direction and 
at most $n$ in the $y$-direction.
We use $k=1$ and $k=2$ for the MHD simulations with $M_A=2.5$ and $M_A=5$, and 
use $k=3$ for the hydrodynamic simulation ($M_A=+\infty$).

\begin{figure}[ht!]
 \caption{\textbf{Example 3:} Contour of vorticity $\omega_h:=\nabla_h\times \bld u_h$ 
 for the test with $M_A = 2.5$.
 Left: $\bar t = 3$, right: $\bar t = 6$. 
 Top: $\mathbf{RT}_1$ scheme, Bottom: $\mathbf{RT}_2$ scheme.
Grid size: $256\times 256$.
}
 \label{fig:kh}
  \includegraphics[width=.48\textwidth]{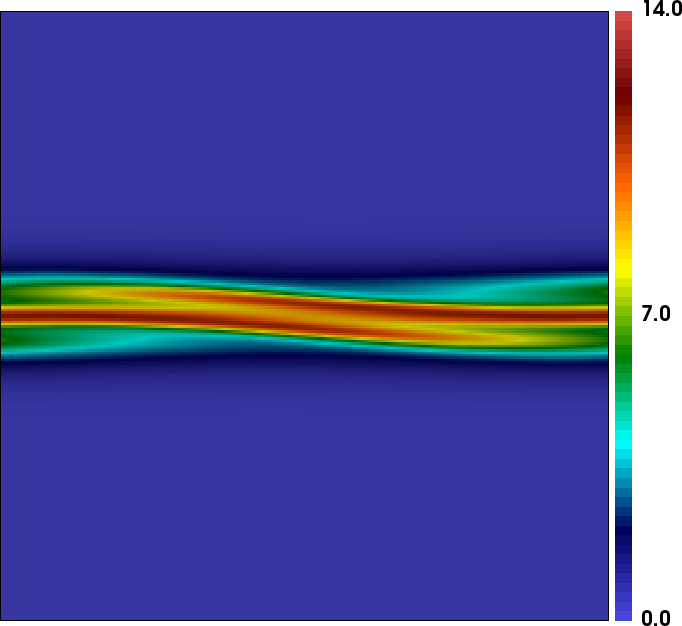}
  \includegraphics[width=.48\textwidth]{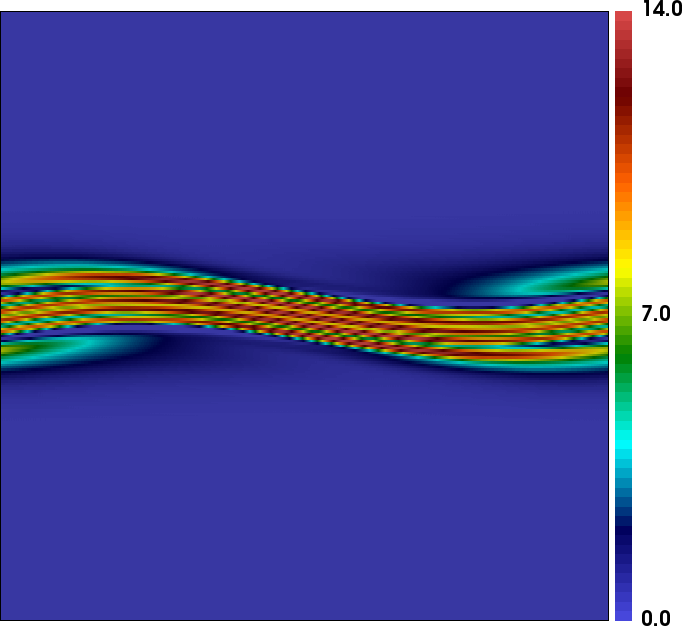}
 \\[1.0ex]
  \includegraphics[width=.48\textwidth]{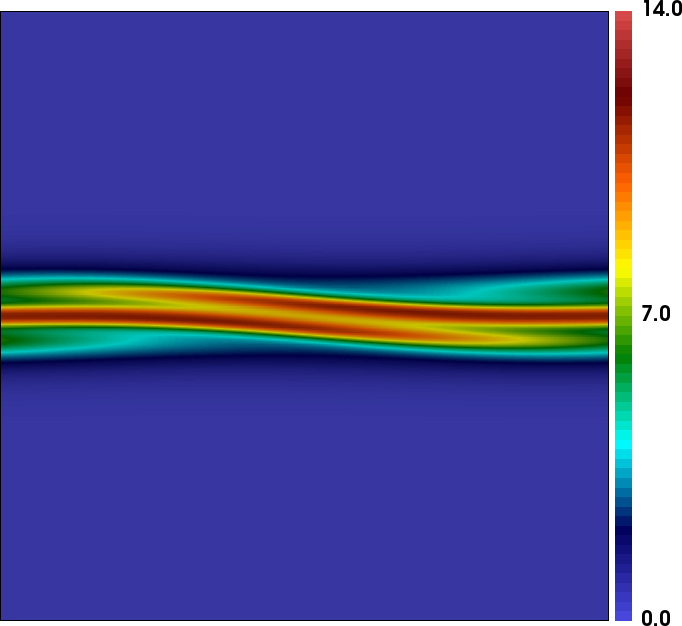}
  \includegraphics[width=.48\textwidth]{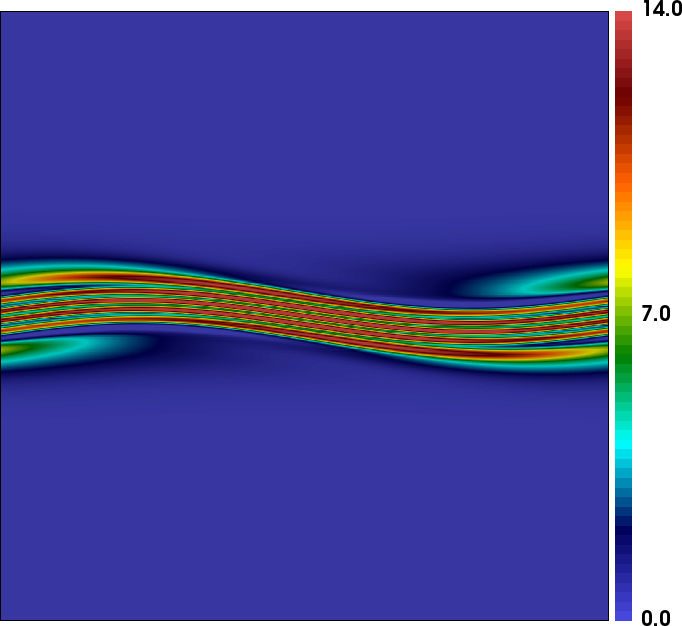}
\end{figure}

\begin{figure}[ht!]
 \caption{\textbf{Example 3:} Contour of vorticity $\omega_h:=\nabla_h\times \bld u_h$ 
 for the test with $M_A = 5$.
 Left: $\bar t = 3$, right: $\bar t = 6$. 
  Top: $\mathbf{RT}_1$ scheme, Bottom: $\mathbf{RT}_2$ scheme.
Grid size: $256\times 256$.
}
 \label{fig:kh2}
  \includegraphics[width=.48\textwidth]{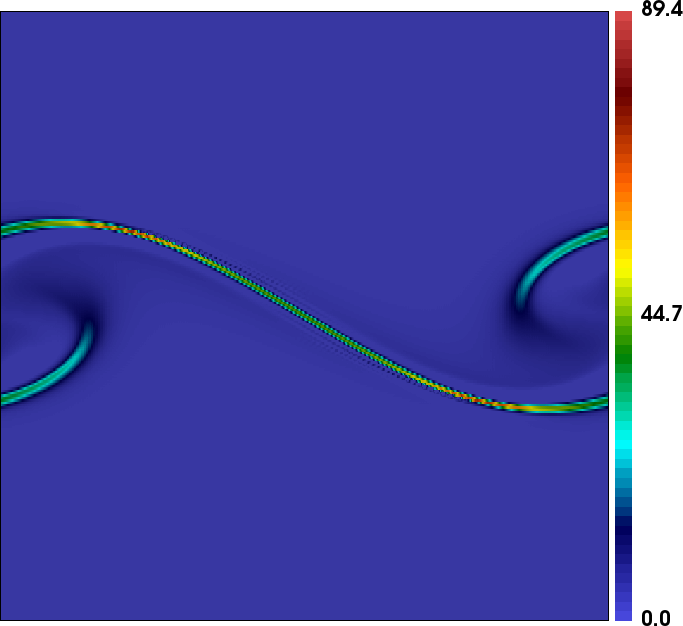}
  \includegraphics[width=.48\textwidth]{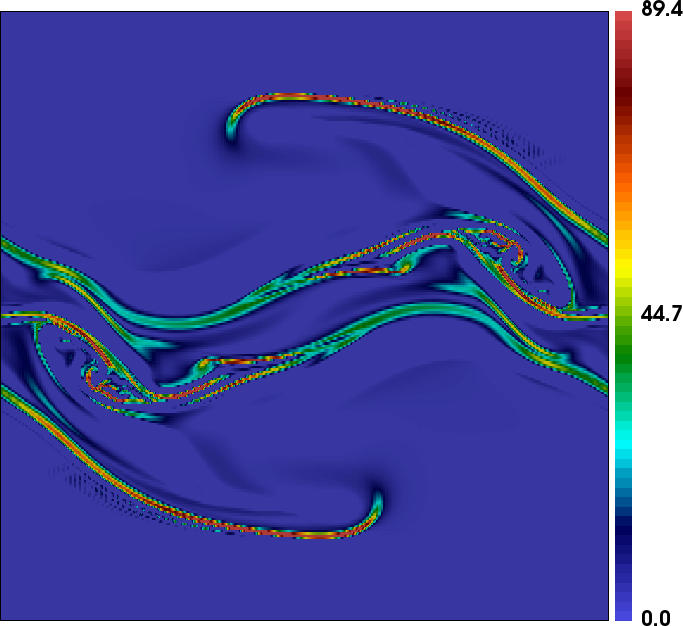}
 \\[1.0ex]
  \includegraphics[width=.48\textwidth]{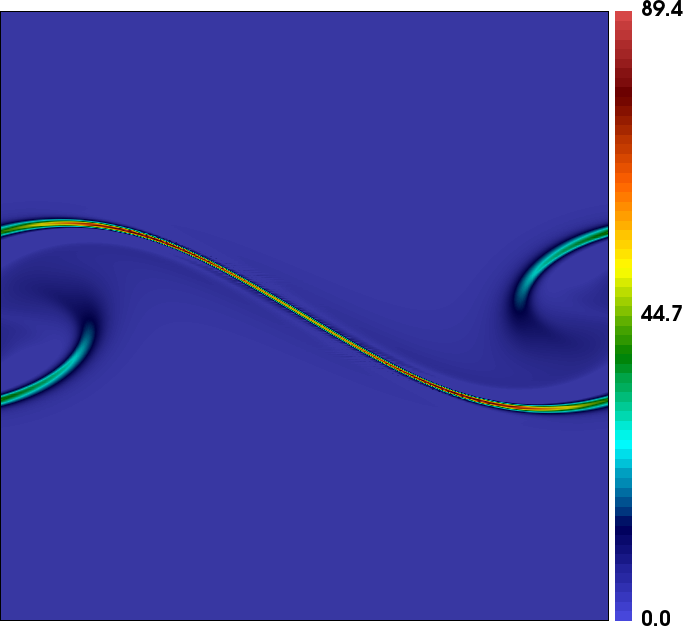}
  \includegraphics[width=.48\textwidth]{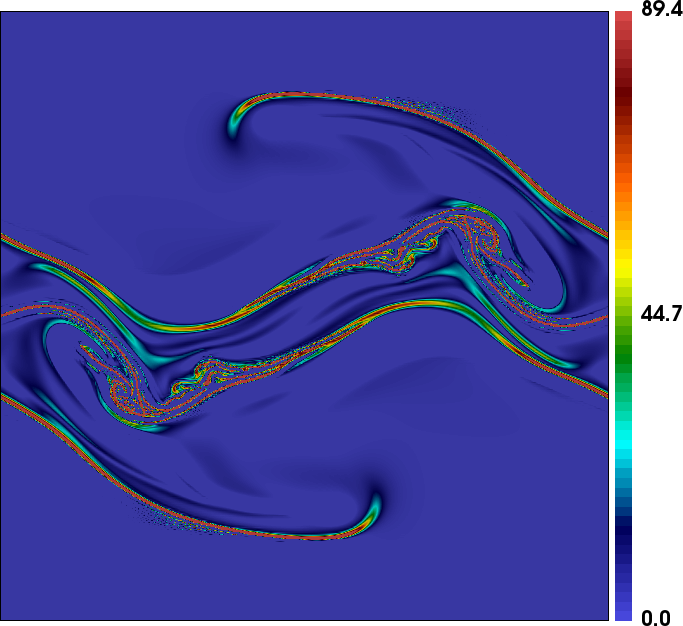}
\end{figure}

\begin{figure}[ht!]
 \caption{\textbf{Example 3:} Contour of vorticity $\omega_h:=\nabla_h\times \bld u_h$ 
 for the test with $M_A = +\infty$ (pure hydrodynamical simulation).
$\mathbf{RT}_3$ scheme. Grid size: $256\times 256$.
}
 \label{fig:kh3}
  \includegraphics[width=.48\textwidth]{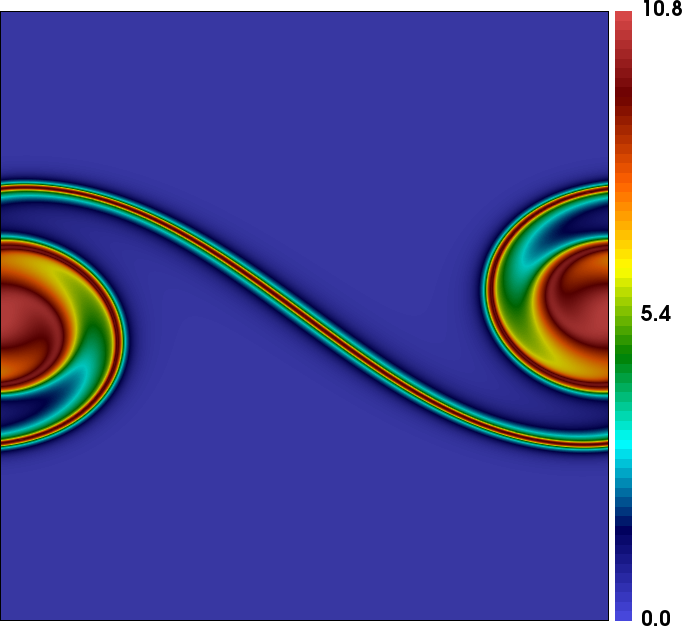}
  \includegraphics[width=.48\textwidth]{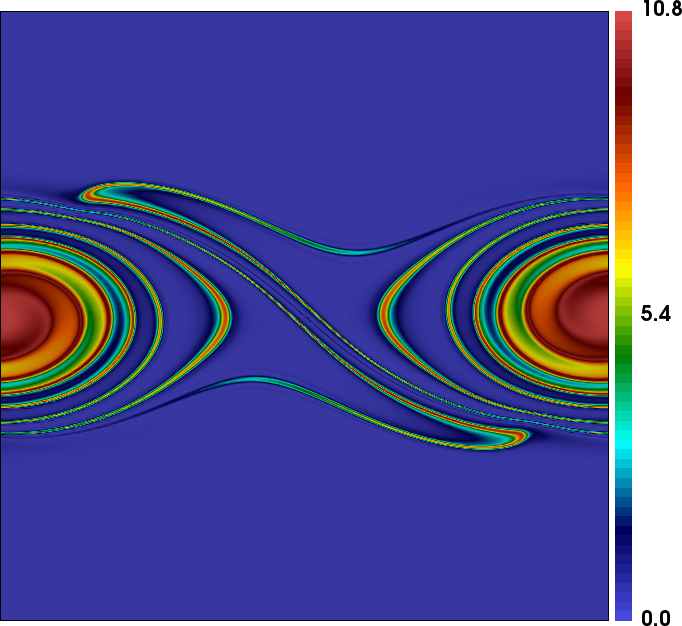}
\end{figure}

Contour of vorticity $\omega_h:=\nabla_h\times \bld u_h$ 
for the cases $M_A=2.5$, $M_A=5$, and $M_A=+\infty$ are 
shown in Fig.~\ref{fig:kh}--\ref{fig:kh3}, respectively.
Comparing with results for the hydrodynamic case in Fig.~\ref{fig:kh3},
we observe that 
in the strong magnetic field case (Fig.~\ref{fig:kh}),
the vortex formation is completely suppressed.
However,  in the weak magnetic field case (Fig.~\ref{fig:kh2}), the 
vortex is initially been developed (left of Fig.~\ref{fig:kh2}),
then destroyed (due to locally strong magnetic field), which, 
in turn, produces a sequence of intermediate vortices.
The flow is significantly more complex for the weak magnetic field case in Fig.~\ref{fig:kh2} than
that for the strong magnetic field case in Fig.\ref{fig:kh}.
These observations are qualitatively in agreement 
with those in \cite{Frank96} for the compressible MHD simulations.
Moreover, it is clear to observe from Fig.~\ref{fig:kh} and Fig.~\ref{fig:kh2}
the resolution improvement from
the second order $\mathbf{RT}_1$ scheme to 
the third order $\mathbf{RT}_2$ scheme.

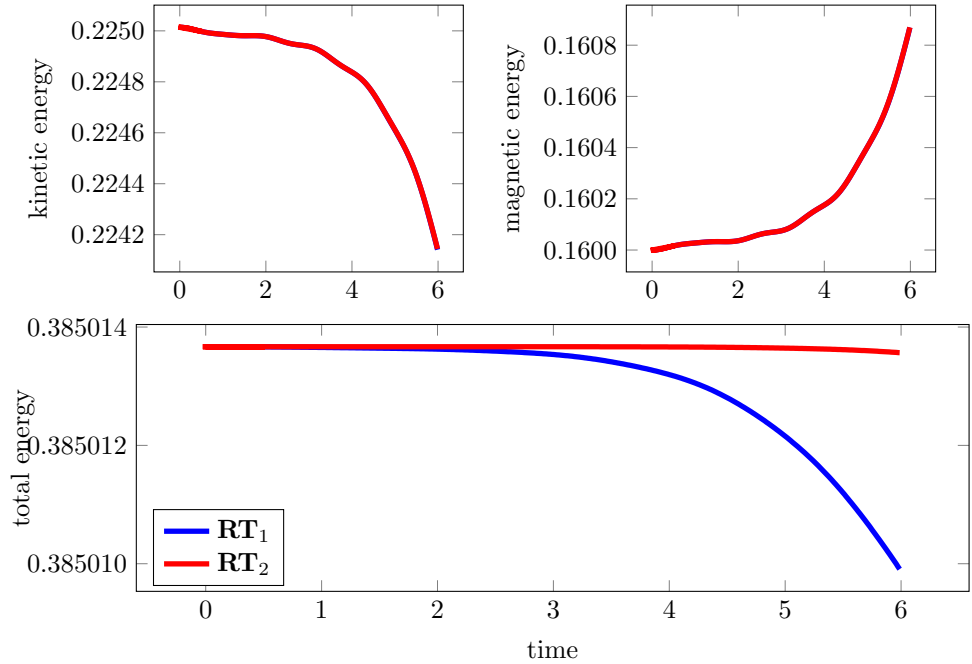
\begin{figure}[ht!]
 \caption{\textbf{Example 3:} 
the time history of kinetic, magnetic, and total energy 
for the test with $M_A=2.5$. 
 }
 \vspace{1ex}
 \label{fig:kh-e1}
 \begin{tikzpicture} 
 \begin{axis}[
 	width=0.45\textwidth,
 	height=0.25\textheight,
    	ylabel={kinetic energy},
	y label style={at={(axis description cs:0.0,.5)},anchor=south},
	  max space between ticks=50pt,
  y tick label style={
    /pgf/number format/.cd,
    fixed,
    fixed zerofill,
    precision=4
  },
    	every axis plot/.append style={line width=2pt, smooth},
    	no markers,   	
    	legend style={at={(0.02,0.02)},anchor=south west}
 	]
\addplot table[x index=0, y index=1]{mkhp1xf2t.csv};
\addplot table[x index=0, y index=1]{mkhp2xf2t.csv};
 \end{axis}
 \end{tikzpicture}
 \hspace{1ex}
  \begin{tikzpicture} 
 \begin{axis}[
 	width=0.45\textwidth,
 	height=0.25\textheight,
    	ylabel={magnetic energy},
    		  max space between ticks=50pt,
  y tick label style={
    /pgf/number format/.cd,
    fixed,
    fixed zerofill,
    precision=4
  },
      y label style={at={(axis description cs:0.0,.5)},anchor=south},
    	every axis plot/.append style={line width=2pt, smooth},
    	no markers,   	
    	legend style={at={(0.02,0.02)},anchor=south west}
 	]
\addplot table[x index=0, y index=2]{mkhp1xf2t.csv};
\addplot table[x index=0, y index=2]{mkhp2xf2t.csv};
 \end{axis}
 \end{tikzpicture}
  \hspace{.5ex}
 \begin{tikzpicture} 
\begin{axis}[
	width=1\textwidth,
	height=0.25\textheight,
	xlabel={time},
	ylabel={total energy},
	y tick label style={/pgf/number format/.cd,/pgf/number format/precision=6,
	/pgf/number format/zerofill=true},
	      y label style={at={(axis description cs:0.0,.5)},anchor=south},
   	every axis plot/.append style={line width=2pt, smooth},
   	no markers,   	
   	legend style={at={(0.02,0.02)},anchor=south west}
	]
\addplot table[x index=0, y index=3]{mkhp1xf2t.csv};
\addplot table[x index=0, y index=3]{mkhp2xf2t.csv};
\addlegendentry{$\mathbf{RT}_1$}
\addlegendentry{$\mathbf{RT}_2$}
\end{axis}
\end{tikzpicture}
\end{figure}

\begin{figure}[ht!]
 \caption{\textbf{Example 3:} 
the time history of kinetic, magnetic, and total energy for the test with $M_A=5$. 
 }
 \vspace{1ex}
 \label{fig:kh-e2}
 \begin{tikzpicture} 
 \begin{axis}[
 	width=0.45\textwidth,
 	height=0.25\textheight,
    	ylabel={kinetic energy},
	y label style={at={(axis description cs:0.0,.5)},anchor=south},
	  max space between ticks=50pt,
  y tick label style={
    /pgf/number format/.cd,
    fixed,
    fixed zerofill,
    precision=2
  },
    	every axis plot/.append style={line width=2pt, smooth},
    	no markers,   	
    	legend style={at={(0.02,0.02)},anchor=south west}
 	]
\addplot table[x index=0, y index=1]{mkhp1yf2t.csv};
\addplot table[x index=0, y index=1]{mkhp2yf2t.csv};
 \end{axis}
 \end{tikzpicture}
 \hspace{1ex}
  \begin{tikzpicture} 
 \begin{axis}[
 	width=0.45\textwidth,
 	height=0.25\textheight,
    	ylabel={magnetic energy},
    		  max space between ticks=50pt,
  y tick label style={
    /pgf/number format/.cd,
    fixed,
    fixed zerofill,
    precision=2
  },
      y label style={at={(axis description cs:0.0,.5)},anchor=south},
    	every axis plot/.append style={line width=2pt, smooth},
    	no markers,   	
    	legend style={at={(0.02,0.02)},anchor=south west}
 	]
\addplot table[x index=0, y index=2]{mkhp1yf2t.csv};
\addplot table[x index=0, y index=2]{mkhp2yf2t.csv};
 \end{axis}
 \end{tikzpicture}
  \hspace{.5ex}
 \begin{tikzpicture} 
\begin{axis}[
	width=1\textwidth,
	height=0.25\textheight,
	xlabel={time},
	ylabel={total energy},
	y tick label style={/pgf/number format/.cd,/pgf/number format/precision=3,
	/pgf/number format/zerofill=true},
	      y label style={at={(axis description cs:0.0,.5)},anchor=south},
   	every axis plot/.append style={line width=2pt, smooth},
   	no markers,   	
   	legend style={at={(0.02,0.02)},anchor=south west}
	]
\addplot table[x index=0, y index=3]{mkhp1yf2t.csv};
\addplot table[x index=0, y index=3]{mkhp2yf2t.csv};
\addlegendentry{$\mathbf{RT}_1$}
\addlegendentry{$\mathbf{RT}_2$}
\end{axis}
\end{tikzpicture}
\end{figure}
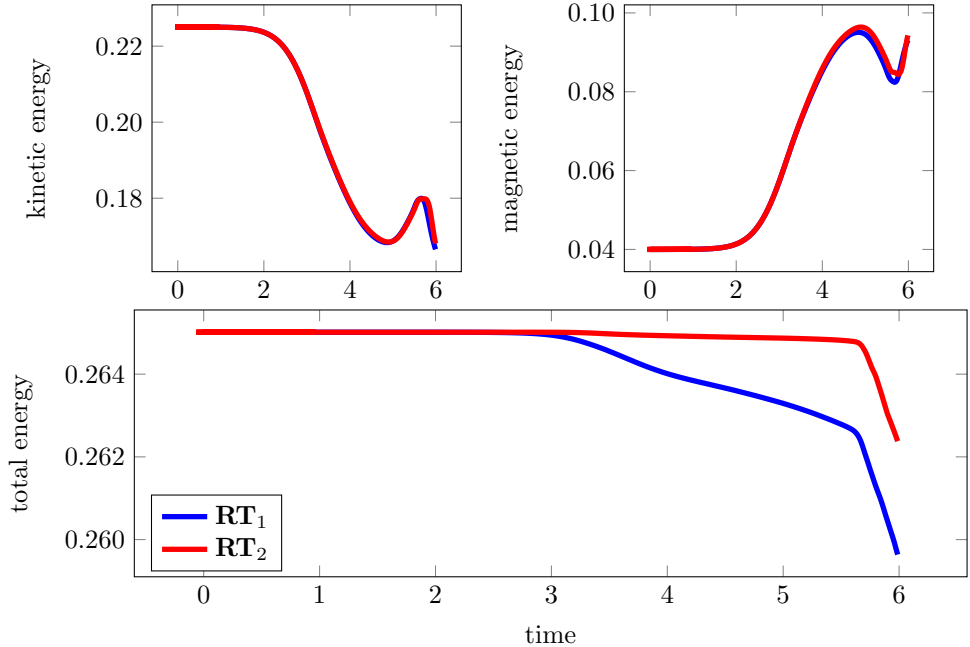

Then, in Fig.~\ref{fig:kh-e1} and Fig.~\ref{fig:kh-e2}, we plot the time evolution of
kinetic, magnetic, and total energy for the case with $M_A=2.5$ and $M_A=5$, respectively.
Looking at the evolution of the total energy, we observe less 
numerical dissipation for the third order $\mathbf{RT}_2$ scheme over the
second order $\mathbf{RT}_1$ scheme, just as expected.
We also observe from Fig.~\ref{fig:kh-e1} that there is no significant energy 
transformation between kinetic and magnetic energy for the strong magnetic field case $M_A=2.5$.
On the other hand, we see from Fig.~\ref{fig:kh-e2} that
the energy transformation for the weak magnetic field case 
$M_A=5$ is quite more complex.
In particular, the kinetic energy evolves through four phases:
it stays at around the same level till scaled time $\bar t\approx 2$, 
then decays till scaled time $\bar t  \approx 5$, then 
increases till scaled time $\bar t \approx 5.5$, and 
then decays till the final scaled time $\bar t = 6$.

\section{Conclusion}
We have presented a DG scheme for the incompressible MHD flow based on a velocity-magnetic field formulation.
Both velocity and magnetic field are discretized using global divergence-free finite elements.
Highlights of the scheme includes global and local conservation properties, high-order
accuracy, energy-stability, and pressure-robustness.
The semidiscrete DG scheme is coupled with explicit time stepping, where in each time step, 
a hybrid-mixed Poisson solver is used to circumvent the (global) mass matrix inversion.
The extension of the current scheme to compressible flow which 
use a divergence-conforming momentum field approximation 
consists of our ongoing work, with the goal of obtaining an asymptotic preserving 
scheme in the low Mach number limit.

{{\bf Acnowledgements}. 
The author would like to thank Prof. Chi-Wang Shu for 
suggesting to work on the problem, and for many helpful discussions concerning the subject.
Part of this research was conducted using computational
resources and services at the Center for Computation and Visualization, Brown University.
}
\bibliographystyle{siam}

\end{document}